\documentclass[11pt]{amsart}
\usepackage{amssymb}
\usepackage{epsfig}
\usepackage{url}
\theoremstyle{plain}

\newtheorem{thm}{Theorem}[section]
\newtheorem{cor}[thm]{Corollary}
\newtheorem{lem}[thm]{Lemma}

\newtheorem{rem}[thm]{Remark}
\newtheorem{ques}[thm]{Question}

\numberwithin{equation}{section}

\def\cal{\mathcal}
\def\bbb{\mathbb}
\def\op{\operatorname}
\renewcommand{\phi}{\varphi}
\newcommand{\R}{\bbb{R}}

\newcommand{\Z}{\bbb{Z}}
\newcommand{\Q}{\bbb{Q}}
\newcommand{\F}{\bbb{F}}

\newcommand{\C}{\bbb{C}}

\begin{document}
\title{Rational points on certain hyperelliptic curves over finite fields}
\author{Maciej Ulas}
\thanks{The author is a scholarship holder of the Stanis\l aw Estreicher Fund.}


\date{}

\begin{abstract}
Let $K$ be a field, $a,\;b\in K$ and $ab\neq 0$. Let us consider
the polynomials $g_{1}(x)=x^n+ax+b,\;g_{2}(x)=x^n+ax^2+bx$, where
$n$ is a fixed positive integer. In this paper we show that for
each $k\geq 2$ the hypersurface given by the equation
\begin{equation*}
S_{k}^{i}:\;u^2=\prod_{j=1}^{k}g_{i}(x_{j}),\quad i=1,\;2.
\end{equation*}
contains a rational curve. Using the above and Woestijne's recent
results \cite{Woe} we show how one can construct a rational point
different from the point at infinity on the curves
$C_{i}:y^2=g_{i}(x),\;(i=1,\;2)$ defined over a finite field, in
polynomial time.
\end{abstract}

\maketitle
    \begin{quote}
      \hskip 3.5cm         {\it Dedicated to the memory of Andrzej M\c{a}kowski}
    \end{quote}
\bigskip

\section{Introduction}

R. Schoof in \cite{Sch} showed how one can count rational points
on the elliptic curve $E: y^2=x^3+ax+b$ defined over finite field
$\F_{p}$, where $p>3$ is a prime, in polynomial time.
Surprisingly, this algorithm allows to compute the order of the
group $E(\F_{p})$ without providing any point (different from the
point at infinity) on the curve $E$ explicite. In this paper a
problem was posed to construct an algorithm which would allow to
determine a rational point $P\in E(\F_{p})\setminus\{\cal{O}\}$ in
polynomial time.

To author's best knowledge, the first work concerning this problem
appeared in 2004. A. Schinzel and M. Ska\l ba showed in
\cite{SchSka} how to determine efficiently a rational point on the
curves of the form $y^2=x^n+a$, where $n=3,\;4,\;a\in\F_{q}$ and
$q=p^{m}$. In case $n=3$, the authors give explicit elements
$y_{1},\;y_{2},\;y_{3},\;y_{4}\in\F_{q}$ with the property that
for at least one $i\leq 4$ the equation $y_{i}^2=x^3+a$ has a
solution in the field $\F_{q}$. In case $n=4$, the authors give a
construction of elements $y_{1},\;y_{2},\;y_{3}\in\F_{q}$ with the
property that for at least one $i\leq 3$, the equation
$y_{i}^2=x^4+a$ has a solution in the field $\F_{q}$.

A solution of the subproblem of finding the $x$-coordinate of a
rational point on the elliptic curve
\begin{equation*}
E:\;y^2=x^3+ax+b=:f(x),
\end{equation*}
 in case when $a,\;b\in \F_{q},\;a\neq 0,\;q=p^{m}$ and
$\F_{q}$ is a finite field with $p>3$, was provided in
\cite{Ska}(for many applications the $y$-coordinate is not
needed). The key element of the proof was a construction of
non-constant rational functions $x_{1},\;x_{2},\;x_{3},\;u\in
K(t)$ which satisfy the equation
\begin{equation}\label{R1}
u^2=f(x_{1})f(x_{2})f(x_{3}).
\end{equation}
We know that multiplicative group $\F_{q}^{\ast}$ is cyclic. This
fact plus the obtained parametric solution prove that for at least
one $i\leq 3$, the element $f(x_{i})$ is a square in $\F_{q}$. If
now $q=p$ or $q=p^m$ and an element $v\in \mathbb{F}_{q}\setminus
\mathbb{F}_{q}^{2}$ is given, then using Schoof's and
Tonneli-Shanks' algorithm (given in \cite{Sha}) respectively, we
can calculate a square root from $f(x_{i})$ in polynomial time.

In his PhD dissertation \cite{Woe}, Ch. van de Woestijne showed
how in polynomial time, for given
$b_{0},\;b_{1},\;\ldots,\;b_{n}\in \F_{q}^{\ast}$ we can find
integers $i,\;j$ such that $0\leq i<j \leq n$, and an element
$b\in\F_{q}^{\ast}$ such that $b_{i}/b_{j}=b^{n}$. We should pay
attention that for calculating this $n$-th root, it is not
necessary to have the element $v\in \mathbb{F}_{q}\setminus
\mathbb{F}_{q}^{n}$ and that the algorithm which computes this
root is deterministic. It is easy to see that having elements
$x_{1},\;x_{2},\;x_{3}$ fulfilling the identities (\ref{R1}) for a
certain $u\in\F_{q}$ and using Woestijne's result we can find a
rational point on the curve $y^2=f(x)$ in polynomial time. This
idea was used in \cite{ShaWoe}. The authors constructed a rational
curve (different from the one constructed in \cite{Ska}) on the
hypersurface $u^2=f(x_{1})f(x_{2})f(x_{3})$, where $f$ is a given
polynomial of degree three. However, we should note that in order
to obtain an explicit form of the curve, it is necessary to solve
the equation $\alpha x^2+ \beta y^2=\gamma$ in $\F_{q}$ for
certain $\alpha,\;\beta,\;\gamma$ (this can be done in
deterministic polynomial time, but of course this lenghtens the
time needed to compute rational point on the curve $E$). The
authors also showed how to construct rational points on elliptic
curves defined over finite fields of characteristic 2 and 3.

A natural question arises here concerning the existence of
rational curves on the hypersurface of the form
\begin{equation}\label{R2}
S_{k}:\;u^2=\prod_{i=1}^{k}g(x_{i}),
\end{equation}
where $g\in\Z[x]$ and $g$ is without multiple roots. It is worth
noting that in this case $S_{k}$ is smooth. It appears that the
problem posed in this form has not been considered so far and
seems to be interesting. Papers \cite{Ska} and \cite{ShaWoe} also
show that in case when $k$ is odd, the ability to construct
rational curves on the hypersurface $S_{k}$ can be useful in
finding rational points on hyperelliptic curves (defined over
finite field) of the form
\begin{equation*}
C:\;y^2=g(x).
\end{equation*}

It is also worth noting that in case when $\op{deg}g=2,3,4$ and
$k$ is an even number, it is easy to find rational curves on
$S_{k}$. Indeed, if $k=2$, then on the surface $S_{2}$ we have a
rational curve $(x_{1},\;x_{2},\;u)=(t,\;t,\;g(t))$. If now
$\op{deg}g=2$, then using a standard procedure we parametrize
rational solutions of equation $u^2=g(t)g(x)$. We act similarly in
case when $\op{deg}g=3$ or $\op{deg}g=4$ with such a difference
that this time we use an algorithm of adding points on curve with
genus one with known rational points. In this way we obtain
infinitely many rational curves on the surface $S_{2}$. As an
immediate consequence of the above reasoning, we obtain curves on
$S_{k}$ in case when $k>2$ is an even integer.

However, if $\op{deg}g>4$ or $k$ is an odd integer, then the task
seems to be much more difficult and the crucial question arises
whether for a given $g\in\Z[x]$ there is a $k$ such that $S_{k}$
contains a rational curve?

\bigskip

Let now $a,\;b\in K,\;ab\neq 0$ and let us consider polynomials
\begin{align*}
&g_{1}(x)=x^n+ax+b, \\
&g_{2}(x)=x^n+ax^2+bx,
\end{align*}
where $n$ is a fixed positive integer. In this paper we prove that
if $g=g_{1}$ or $g=g_{2}$, then for each $k\geq 2$ there is a
rational curve on the surface $S_{k}$.

These and Woestijne's results show that the construction of
rational points on the curve $C_{i}: y^2=g_{i}(x)$ can be
performed in polynomial time. Let us also note that if $n$ is
even, then $g_{i}(-b/a)=(b/a)^{n}$ and the point
$P=(-b/a,\;(b/a)^{n/2})$ lies on the curve $y^2=g_{i}(x)$ and is
different from the point at infinity. One can see that in this
case the problem of existence of rational points on $C_{i}$ is
easy. Certainly, it does not provide us with an answer to the
question about constructing a rational curve on hypersurface
$S_{k}$ when $n$ is even.

\bigskip
\section{Rational curves on $ S_{k}^{i}$}

In this section we consider the hypersurface
\begin{equation*}
S_{k}^{i}:\;u^2=\prod_{j=1}^{k}g_{i}(x_{j}),
\end{equation*}
where $i\in\{1,\;2\}$ is fixed. As a direct examination of the
existence of rational curves on the hypersurace $S_{k}^{i}$ is
difficult, let us reduce our problem to the examination of simpler
objects.

\bigskip

 Let $a,\;b,\;c,\;d\in K$ fulfill the condition
\begin{equation*}
(\ast)\quad a\neq 0 \quad\mbox{or}\quad c\neq 0
\quad\mbox{and}\quad b\neq 0 \quad\mbox{or}\quad d\neq 0.
\end{equation*}
Let now $m,\;n$ be fixed positive integers and let us consider the
surfaces
\begin{align*}
&S^{1}:\;g_{1}(x)z^m=y^n+cy+d,\\
&S^{2}:\;g_{2}(x)z^m=y^n+cy^2+dy.
\end{align*}

We will prove that a rational curve lies on each of these
surfaces. Using constructed curves we will show how to construct
curves on $S_{2}^{i}$ and $S_{3}^{i}$. Because each positive
integer $\geq 2$ is of the form $2k+3l$, then as an immediate
consequence we obtain the existence of rational curves on
hypersurface $S_{k}^{i}$ for each $k\geq 2$.

\bigskip

We start with the following

\begin{lem}\label{lem1}
Let $n,\;m\in\mathbb{N}_{+}$ and $a,\;b,\;c,\;d\in K$ fulfill
condition $(\ast)$. Then on each of the surfaces $S^{1},\;S^{2}$
there is a rational curve.
\end{lem}

\begin{proof} Let $F_{1}(x,y,z):=g_{1}(x)z^m-(y^n+cy+d)$ and
$F_{2}(x,y,z):=g_{2}(x)z^m-(y^n+cy^2+dy)$. For the proof, let us
put $x=T,\;y=t^mT,\;z=t^n$. It is easy to see that for $x,\;y,\;z$
defined in this way, the equation $F_{1}(T,t^mT,t^n)=0$ has the
root
\begin{equation*}
T=-\frac{bt^{mn}-d}{at^{mn}-ct^{m}},
\end{equation*}
which gives us a parametric curve $L_{1}$ on the surface $S^{1}$
given by the equations:
\begin{equation*}
L_{1}:\quad x(t)=-\frac{bt^{mn}-d}{at^{mn}-ct^{m}},\quad
y(t)=-\frac{bt^{mn}-d}{at^{m(n-1)}-c},\quad z(t)=t^n.
\end{equation*}

The same method can be applied to find a rational curve on the
surface $S^{2}$. In this case the equation $F_{2}(T,t^mT,t^n)=0$
has two roots, $T=0$ and
\begin{equation*}
T=-\frac{bt^{m(n-1)}-d}{at^{m(n-1)}-ct^{m}}.
\end{equation*}
Rational curve $L_{2}$ on $S^{2}$ is given by the equations
\begin{equation*}
L_{2}:\quad x(t)=-\frac{bt^{m(n-1)}-d}{at^{m(n-1)}-ct^{m}},\quad
y(t)=-\frac{bt^{m(n-1)}-d}{at^{m(n-2)}-c},\quad z(t)=t^n.
\end{equation*}

Note that the condition $(\ast)$ plays a crucial role in our
reasoning in both cases.

\end{proof}

\begin{rem}\label{rem1}
{\rm The surface $S^{1}$ appeared in \cite{Mes} with an additional
assumption $a=c,\;b=d,\;m=2,\;n=3$. In this case the curve $L_{1}$
was used to show that for a given $j\neq 0,\;1728$ there are
infinitely many elliptic curves with $j$-invariant equal to $j$
and Mordell-Weil rank $\geq 2$.

It should be noted that a special case of the surface $S^{1}$,
when $m=2,\;n=3$, was also considered in \cite{KuWa}. In this case
the curve $L_{1}$ was used to show that on the surface $S^{1}$ the
set of rational points is dense in the topology of $\R^{3}$.}
\end{rem}

\bigskip

Using the above lemma we can prove the following

\begin{thm}\label{thm2}
Let $K$ be a field and let us put
$g_{1}(x)=x^n+ax+b,\;g_{2}(x)=x^n+ax^2+bx$, where $a,\;b\in
K,\;ab\neq 0$. Let $t,\;u$ be variables.
\begin{enumerate}
\item  If $n\geq 3$ is a fixed positive integer and let us put
\begin{equation*}
X_{1}(t)=-\frac{b}{a}\frac{t^{2n}-1}{t^{2n}-t^2},\quad
X_{2}(t)=t^2X_{1}(t),\quad U(t)=t^ng_{1}(X_{1}(t)).
\end{equation*}
Then
\begin{equation*}
U(t)^2=g_{1}(X_{1}(t))g_{1}(X_{2}(t)).
\end{equation*}
If now
\begin{equation*}
X_{1}(t)=-\frac{b}{a}\frac{t^{2(n-1)}-1}{t^{2(n-1)}-t^2},\quad
X_{2}(t)=t^2X_{1}(t),\quad U(t)=t^ng_{1}(X_{1}(t)),
\end{equation*}
then
\begin{equation*}
U(t)^2=g_{2}(X_{1}(t))g_{2}(X_{2}(t)).
\end{equation*}
\item If $n$ is an odd integer, let us put
\begin{align*}
X_{1}(t,u)&=u,\\
X_{2}(t,u)&=-\frac{b}{a}\frac{t^{2n}g_{1}(u)^n-1}{g_{1}(u)(t^{2n}g_{1}(u)^{n-1}-t^2)},\\
X_{3}(t,u)&=t^2g_{1}(u)X_{2}(t,u),\\
U(t,u)&=t^ng_{1}(u)^{(n+1)/2}g_{1}(X_{2}(t,u)).
\end{align*}
Then
\begin{equation*}
U(t,u)^2=g_{1}(X_{1}(t,u))g_{1}(X_{2}(t,u))g_{1}(X_{3}(t,u)).
\end{equation*}
If now
\begin{align*}
X_{1}(t,u)&=u,\\
X_{2}(t,u)&=-\frac{b}{a}\frac{t^{2(n-1)}g_{2}(u)^{n-1}-1}{g_{2}(u)(t^{2(n-1)}g_{2}(u)^{n-2}-t^2)},\\
X_{3}(t,u)&=t^2g_{2}(u)X_{2}(t,u),\\
U(t,u)&=t^ng_{2}(u)^{(n+1)/2}g_{2}(X_{2}(t,u)),
\end{align*}
then
\begin{equation*}
U(t,u)^2=g_{2}(X_{1}(t,u))g_{2}(X_{2}(t,u))g_{2}(X_{3}(t,u)).
\end{equation*}
\end{enumerate}
\end{thm}
\begin{proof} We consider the surfaces $S^{1}$ and $S^{2}$ from the Lemma \ref{lem1} with $m=2$.
For the proof of the first part of our theorem, let us make a change of variables
$z=z_{1}/g_{1}(x)$ in the equation of the surface $S^{1}$. By this
change of variables the surface $S^{1}$ is birational with the
surface given by the equation
\begin{equation}\label{R3}
S':\;z_{1}^2=(x^n+ax+b)(y^n+cy+d).
\end{equation}
Putting now $a=c,\;b=d$ and using equations of the curve $L_{1}$
from the proof of Lemma \ref{lem1}, we obtain the statement of our
theorem.

Now we take the equations defining the curve $L_{2}$ from the
proof of Lemma \ref{lem1} and repeat the above reasoning in case
of the surface $S^{2}$. This ends the proof of the first part of
our theorem.

\bigskip

For the proof of the second part of our theorem let us go back to
the surface $S'$ given by the equation (\ref{R3}). If we now put
$c=a/g_{1}(u)^{n-1},\;d=b/g_{1}(u)^{n}$ and perform a change of
variables
\begin{equation}\label{R4}
u=X_{1},\quad x=X_{2},\quad y=\frac{X_{3}}{g_{1}(u)},\quad
z_{1}=U_{1}g_{1}(u)^{-\frac{n+1}{2}},
\end{equation}
then after elementary calculations the equation of the surface
$S'$ is of the form
\begin{equation*}
U_{1}^2=g_{1}(X_{1})g_{1}(X_{2})g_{1}(X_{3}).
\end{equation*}
If now $x,\;y,\;z$ are rational functions defining the curve
$L_{1}$ on the surface $S^{1}$ for
$c=a/g_{1}(u)^{n-1},\;d=b/g_{1}(u)^{n}$, then calculating
$X_{1},\;X_{2},\;X_{3}$ from equations (\ref{R4}), we obtain a
two-parametric solution of the above equation given by the
expressions from the statement of our theorem.

\bigskip

 The proof of the second part of our theorem in case of the polynomial
$g_{2}(x)=x^n+ax^2+bx$ is similar, with one difference: we should
substitute $a/g_{2}(u)^{n-2}$ and $b/g_{2}(u)^{n-1}$ in place of
$c,\;d$ respectively.
\end{proof}

\begin{rem}\label{rem2}
{\rm Let us note that in case when $K$ is a finite field with
$\op{char}K>3,\;a,\;b\in K,\;ab\neq 0$ and we are interested in
construction of a rational point on the elliptic curve
\begin{equation*}
E:\;y^2=x^3+ax+b=:f(x),
\end{equation*}
then our rational curve lying on the hypersurface
$S:\;u^2=f(x_{1})f(x_{2})f(x_{3})$ is much simpler than that
obtained by Ska\l ba. If $x_{i}=X_{i}(t),\;i=1,\;2,\;3$ are the
equations defining the curve on $S$, then if $X_{1}X_{2}X_{3}=N/D$
for certain relatively prime polynomials $N,\;D\in K[t]$, then
$\op{deg}N\leq 26,\;\op{deg}D\leq 25$ in case of parametrization
obtained by Ska\l ba and $\op{deg}N\leq 8,\;\op{deg}D\leq 6$ in
case of our parametrization (with $u\in K$ such that $f(u)\neq 0$)
from Theorem \ref{thm2}. Multiplicative structure of functions
$X_{i}$ is very simple in our case, too, which influences the
speed of calculations.

Our parametrisation has also this advantage over the one obtained
by Shallue and Woestijne that it is not necessary to solve the
equation of the form $\alpha x^2+\beta y^2=\gamma$ in $K$ in order
to obtain it.}
\end{rem}

As in case when $n$ is even we have $g_{i}(-a/b)=((b/a)^{n/2})^2$,
then from the above theorem we obtain

\begin{cor}\label{cor3}
Let $K$ be a field, $a,\;b\in K,\;ab\neq 0$. Then for each
positive integer $k\geq 2$ there is a rational curve on the
hypersurface
\begin{equation*}
S_{k}^{i}:\;u^2=\prod_{j=1}^{k}g_{i}(x_{j}),\quad i=1,\;2.
\end{equation*}

\end{cor}

Because the case $k=3$ and $K=\F_{q},\;q=p^{m}$ is especially
interesting for us, we have to decide about the assumptions
permitting to calculate the values of function $X_{i}(u,t)$ for
$i=1,\;2,\;3$ from the second part of Theorem \ref{thm2}. We can
limit our considerations to examining the case of polynomial
$g_{1}$, of odd degree $n$. In case of the polynomial $g_{2}$ the
reasoning will be similar.

Firstly, let us note that the functions $X_{i}(t,u)$ for $i=1,2,3$
are non-constant. Moreover, functions $X_{2}$ and $X_{3}$ have the
same denominator which equals
\begin{equation*}
D(t,u)=g_{1}(u)t^{2}\frac{(t^2g_{1}(u))^{n-1}-1}{t^2g_{1}(u)-1}.
\end{equation*}

We know that for each $v\in\F_{q}$ we have $v^{p^{m}-1}=1$. There
are $p^m-n$ elements $u\in\F_{q}$ for which $g_{1}(u)\neq 0$. If
we fix such an element now, then because
$\op{deg}_{t}D(t,u)=2(n-1)$ and $t^2\mid D(t,u)$, there are at
least $p^m-2(n-1)+1$ elements $t\in \F_{q}$ for which $D(t,u)\neq
0$. Thus we can see that there are at least
$(p^m-n)(p^m-2(n-1)+1)$ elements in $\F_{q}\times \F_{q}$ for
which $D(t,u)\neq 0$. From this observation we see that if
$p>2(n-1)-1$ then we can find $t,\;u\in \F_{q}$ such that for at
least one $j\in\{1,\;2,\;3\}$ the $g_{1}(X_{j}(t,u))$ is a square.

 \bigskip
\section{Some remarks and questions}

\bigskip

Let us define the set $T$ which contains pairs $(t,u)\in
\F_{q}\times \F_{q}$ for which we can compute
$X_{i}(t,u),\;i=1,2,3$ from the preceding section. Then, we can
define the map $\Phi$ from $T$ to the curve $C:y^2=g_{1}(x)$ in
the following way
\begin{equation*}
\Phi(t,u)=(X_{j}(t,u),\sqrt{g_{1}(X_{j}(t,u))}\;),
\end{equation*}
where square root is taken in $\F_{q}$ and
$j=\op{min}\{i:g_{1}(X_{i}(t,u))\; \mbox{is a square}\}$. Note
that there are at most $2q$ rational points on $C$ over $\F_{q}$
but $T$, as we have proved, contains at least $(q-n)(q-2(n-1)+1)$
elements. This suggest the following

\begin{ques}\label{ques1}
Is the map $\Phi:T\ni (t,u)\longmapsto \Phi(t,u)\in C$ a
surjective map?
\end{ques}

Another question which comes to mind is the following.

\begin{ques}\label{ques2}
Let us fix a polynomial $g\in\Z[x]$ without multiple roots. Is
there an integer $k\geq 2$ such that on the hypersurface
\begin{equation*}
 S_{k}:\;u^2=\prod_{j=1}^{k}g(x_{j})
\end{equation*}
there are infinitely many rational points with $u\neq 0$? In this
question we are interested in non-trivial points on $S_{k}$, i. e.
such points $(x_{1},\;\ldots,\;x_{k},\;u)$ that $g(x_{i})\neq
g(x_{j})$ for $i\neq j$.
\end{ques}

It would be also interesting to know whether if we treat
hypersurface $S_{k}$ over $\C$ (instead $\Q$) then there are
rational curves on $S_{k}$.

It seems that the following question is much more difficult.

\begin{ques}\label{ques3}
Let us fix a polynomial $g\in\Z[x]$ without multiple roots and a
positive integer $k\geq 2$. Is there a non-trivial rational point
with $u\neq 0$ on the hypersurface $S_{k}$?
\end{ques}

If our fixed $k$ in Question \ref{ques3} is odd we should also
assume that for each $p\in \mathbb{P}\cup\{\infty\}$ the curve
$y^2=g(x)$ has a point over $\Q_{p}$ (as usual $\Q_{\infty}=\R)$.
It is clear that the assumption concerning local solubility is
necessary. For example, consider the polynomial $g(x)=3-x^2$.
There are no $\Q_{3}$-rational points on the curve $y^2=g(x)$ and
this immediately implies that there are no $\Q_{3}$-rational
points on $S_{k}$.

\bigskip
Let us note that if the polynomial $g$ fulfills the condition
$x^ng(1/x)=g(x)$ (this is so called reciprocal polynomial), then
we have a rational curve
$x_{1}=t^2,\;x_{2}=1/t^2,\;u=t^{n}g(1/t^2)$ on the surface
$S_{2}$. As an immediate consequence we conclude that if $k$ is
even, then there is a rational curve on $S_{k}$. Additionally, if
the degree of $g$ is odd, then on the surface $S_{3}$ we have a
rational curve given by equations
\begin{equation*}
x_{1}=t,\;x_{2}=g(t),\;x_{3}=\frac{1}{g(t)},\;u=g(t)^{\frac{n+1}{2}}g\Big(\frac{1}{g(t)}\Big),
\end{equation*}
and immediately we have that for each $k\geq 2$ there is a
rational curve on $S_{k}$.

Let us also note that if $g(x)=x^4+1$, then on the hypersurface
$S_{3}$ we have a rational curve with $x_{i}=x_{i}(t),\;i=1,2,3$
given by
\begin{equation*}
x_{1}=\frac{2t+1}{3t^2+3t+1},\;x_{2}=\frac{3t^2+2t}{3t^2+3t+1},\;x_{3}=\frac{3t^2+4t+1}{3t^2+3t+1}
\end{equation*}

\bigskip

It would be very interesting to construct other families of
polynomials of such a property that for each $k\geq 2$ there are
rational curves (or infinitely many non-trivial rational points)
on $S_{k}$.

\bigskip
\noindent {\bf Acknowledgments.} I would like to thank the
anonymous referee for his valuable comments and Professors A.
Schinzel and K. Rusek for remarks improving the presentation. I am
also grateful to Dr M. Ska\l ba who acquainted me with the
Woestijne paper \cite{Woe}.

\bigskip

\hskip 4.5cm       Maciej Ulas

 \hskip 4.5cm       Jagiellonian University

 \hskip 4.5cm       Institute of Mathematics

 \hskip 4.5cm       Reymonta 4

 \hskip 4.5cm       30 - 059 Krak\'{o}w, Poland

 \hskip 4.5cm      e-mail:\;{\tt Maciej.Ulas@im.uj.edu.pl}


\bigskip

\begin{thebibliography}{100}

\bibitem{KuWa}
M. Kuwata, L. Wang, {\it Topology of rational points on isotrivial
elliptic surfaces}, Int. Math. Research Notices, {\bf 4} (1993),
113-123.

\bibitem{Mes} J. F. Mestre, {\it Rang de courbes elliptiques
d'invariant donn\'{e}}, C. R. Acad. Sci. Paris S\'{e}r. I Math.
{\bf 314} (1992), 919-922.

\bibitem{SchSka}
A Schinzel, M. Ska\l ba, {\it On Equations $y^2=x^n+k$ in a Finite
Field,} Bull. Polish Acad. Sci. Math. {\bf 52} (2004), 223-226.

\bibitem{Sch} R. Schoof, {\it Elliptic curves over finite fields
and the computation of square roots $\bmod\;p$}, Math. Comp. {\bf
44} (170) (1985) 483-494.

\bibitem{ShaWoe}
A. Shallue, Ch. van de Woestijne, {\it Construction of Rational
Points on Elliptic Curves over Finite Fields}, ANTS VII, Lecture
Notes in Comput. Sci., Berlin (2006) 510-524.

\bibitem{Sha}
D. Shanks, {\it Five number theoretic algorithms,} Congr. Numer.
{\bf 7} (1972), 51-70.

\bibitem{Ska}
M. Ska\l ba, {\it Points on elliptic curves over finite fields},
Acta Arith. {\bf 117} (2005), 293-301.

\bibitem{Woe} Ch. van de Woestijne, {\it Deterministic equation solving over finite fields}, PhD thesis, Universiteit Leiden (2006).


\end{thebibliography}
 \end{document}